\pdfoutput=1
\RequirePackage{ifpdf}
\ifpdf %We are running pdfTeX in pdf mode
\documentclass[pdftex]{sigma}
\else
\documentclass{sigma}
\fi

\begin{document}

\newcommand{\arXivNumber}{1410.1232}

\allowdisplaybreaks

\renewcommand{\thefootnote}{$\star$}

\renewcommand{\PaperNumber}{044}

\FirstPageHeading

\ShortArticleName{Time and Band Limiting for Matrix Valued Functions, an Example}

\ArticleName{Time and Band Limiting for Matrix Valued Functions, an Example\footnote{This paper is a~contribution to the Special Issue on Exact Solvability and Symmetry Avatars
in honour of Luc Vinet.
The full collection is available at
\href{http://www.emis.de/journals/SIGMA/ESSA2014.html}{http://www.emis.de/journals/SIGMA/ESSA2014.html}}}

\Author{F.~Alberto GR\"UNBAUM~$^\dag$, In\'es PACHARONI~$^\ddag$ and Ignacio Nahuel ZURRI\'AN~$^\ddag$}

\AuthorNameForHeading{F.A.~Gr\"unbaum, I.~Pacharoni and I.N.~Zurri\'an}

\Address{$^\dag$~Department of Mathematics, University of California, Berkeley~94705, USA}
\EmailD{\href{mailto:grunbaum@math.berkeley.edu}{grunbaum@math.berkeley.edu}}

\Address{$^\ddag$~CIEM-FaMAF, Universidad Nacional de C\'ordoba, C\'ordoba~5000, Argentina}
\EmailD{\href{mailto:pacharon@famaf.unc.edu.ar}{pacharon@famaf.unc.edu.ar},
\href{mailto:zurrian@famaf.unc.edu.ar}{zurrian@famaf.unc.edu.ar}}

\ArticleDates{Received February 11, 2015, in f\/inal form May 30, 2015; Published online June 12, 2015}

\Abstract{The main purpose of this paper is to extend to a~situation involving matrix valued orthogonal polynomials and
spherical functions, a~result that traces its origin and its importance to work of Claude Shannon in laying the
mathematical foundations of information theory and to a~remarkable series of papers by D.~Slepian, H.~Landau and
H.~Pollak.
To our knowledge, this is the f\/irst example showing in a~non-commutative setup that a~bispectral property implies that
the corresponding global operator of ``time and band limiting'' admits a~commuting local operator.
This is a~noncommutative analog of the famous prolate spheroidal wave operator.}


\Classification{33C45; 22E45; 33C47}

\renewcommand{\thefootnote}{\arabic{footnote}}
\setcounter{footnote}{0}

\section{Introduction}

In a~ground-breaking paper laying down the mathematical foundations of communication theory, Claude Shannon~\cite{S,S+}
considers a~basic problem in harmonic analysis and signal processing: how to best concentrate a~function both in
physical and frequency space.
This issue was an important part of the work of C.~Shannon for several years after the publication of this paper.
The problem was appeared earlier in several versions and one should at least mention the role of the Heisenberg
inequality in this context: for a~nice and simple proof it~-- due to W.~Pauli~-- see Hermann Weyl's book~\cite{W}.

What is really novel in Shannon and coworkers's look at this problem is the following question: suppose you consider an
unknown signal $f(t)$ of f\/inite duration, i.e., the signal is non-zero only in the interval $[-T,T]$. The data you have
are the values of the Fourier transform $\mathcal Ff(k)$ of~$f$ for values of~$k$ in the interval $[-W,W]$.
What is the best use you can make of this data?

In practice, the values of $\mathcal Ff(k)$ will be corrupted by noise and one is dealing with a~typical situation in
signal processing: recovering an image from partial and noisy data in the presence of some apriori information.

This problem was treated originally by Shannon himself but a~full solution had to wait for the joint work, in dif\/ferent
combinations, of three remarkable workers at Bell labs in the 1960's: David Slepian, Henry Landau and Henry Pollak,
see~\cite{SLP2,SLP3,SLP4,SLP5, SLP1}.

A very good account of this development is a~pair of papers by David Slepian~\cite{S2,S1}.
The f\/irst one is essentially the second Shannon lecture given at the International Symposium on Information Theory in~1975.
The abstract starts with the sentence ``It is easy to argue that real signals must be bandlimited.
It is also easy to argue that they cannot be so''.
The ideas in this paper took their def\/inite form in~\cite{M}.

The second paper is on the occasion of the John von Neumann lecture given at the SIAM 30th anniversary meeting in 1982.
Here is a~quote from the second paper: ``There was a~lot of serendipity here, clearly.
And then our solution, too, seemed to hinge on a~lucky accident-namely that we found a~second-order dif\/ferential
operator that commuted with an integral operator that was at the heart of the problem''.

What these workers found was that instead of looking for the unknown $f(t)$ itself one should consider a~certain
integral operator with discrete spectrum in the open interval $(0,1)$ and a~remarkable ``spectral gap'': about $[4 WT]$
(integer part of $2W \times 2T$) eigenvalues are positive, and all the remaining ones are essentially zero.
They argue that in the presence of noisy data one should try to compute the projection of $f(t)$ on the linear span of
the eigenfunctions with ``large'' eigenvalues.
The ef\/fective computation of these eigenfunctions is made possible by replacing the integral operator by the commuting
dif\/ferential one alluded to by D.~Slepian (both have simple spectrum).
From a~theoretical point of view these eigenfunctions are the same, but using the dif\/ferential operator instead of the
integral one, we have a~manageable numerical problem.
For a~very recent account of several computational issues see~\cite{ORX,JKS} and~\cite{BK14}.
For new areas of applications involving (sometimes) vector-valued quantities on the sphere, see~\cite{JB,PS, SD,SDW}.

We still have to answer the question:
How do you explain the existence of this local commuting operator?

To this day nobody has a~simple explanation for this miracle.
Indeed there has been a~systematic ef\/fort to see if the ``bispectral property'' f\/irst considered in~\cite{DG},
guarantees the commutativity of these two operators, a~global and a~local one.
A~few papers where this question has been taken up, include~\cite{G3,G4,G5,G6,GLP,GY,P1,P2}.

The results in the present paper are a~(f\/irst) natural extension of the work in~\cite{G3,GLP}, where the classical
orthogonal polynomials played a~central role, to a~matrix-valued case involving matrix orthogonal polynomials.
In a~case such as in~\cite{G3,GLP} or in the present paper where physical and frequency space are of a~dif\/ferent nature
(one is continuous and the other one is discrete) one can deal (as explained in~\cite{G3,GLP}) with either an integral
operator or with a~full matrix.
Each one of these global objects will depend on two parameters that play the role of $(T,W)$ in Shannon's case, and one
will be looking for a~commuting local object, i.e., a second-order dif\/ferential operator or a~tridiagonal matrix.
In this paper we will be dealing with a~full matrix and a~(block) tridiagonal one.

Finally a~word about possible applications.
The paper~\cite{GLP} was written with no particular application in mind, but with the expectation that the analysis of
functions def\/ined on the sphere could benef\/it from it.
A~few years later some applications did emerge, see~\cite{SD,SDW} and its references.
One can only hope that the present paper dealing with matrix valued functions def\/ined on spheres will f\/ind a~natural
application in the future.

\section{Preliminaries}

Let $W=W(x)$ be a~weight matrix of size~$R$ in the open interval $(a,b)$.
By this we mean a~complex $R \times R$-matrix valued integrable function~$W$ on the interval $(a, b)$ such that $W (x)$
is positive def\/initive almost everywhere and with f\/inite moments of all orders.
Let $Q_w(x)$ be a~sequence of real valued matrix orthonormal polynomials with respect to the weight $W(x)$.
Consider the following two Hilbert spaces: The space $L^2((a,b), W(t)dt)$, or simply denoted by $L^2(W)$, of all matrix
valued measurable matrix valued functions $f(x)$, $x\in (a,b)$, satisfying $\int_a^b
\operatorname{tr}(f(x)W(x)f^*(x))dx < \infty $ and the space $\ell^2(M_R, {\mathbb N}_0)$ of all real valued
$R\times R$ matrix sequences $(C_w)_{w\in {\mathbb N}_0}$ such that $\sum\limits_{w=0}^\infty \operatorname{tr}
(C_w C_w^*) < \infty$.

The map $\mathcal F\colon \ell^2(M_R,{\mathbb N}_0) \longrightarrow L^2(W)$ given~by
\begin{gather*}
(A_w)_{w=0}^\infty \longmapsto \sum\limits_{w=0}^\infty A_w Q_w(x)
\end{gather*}
is an isometry.
If the polynomials are dense in $L^2(W)$, this map is unitary with the inverse $\mathcal F^{-1}\colon L^2(W)\longrightarrow
\ell^2(M_R,{\mathbb N}_0) $ given~by
\begin{gather*}
f \longmapsto A_w=\int_a^b f(x) W(x) Q^*_w(x) dx.
\end{gather*}

We call our map $\mathcal F$ to remind ourselves of the usual Fourier transform.
Here ${\mathbb N}_0$ takes up the role of ``physical space'' and the interval $(a,b)$ the role of ``frequency space''.
This is, clearly, a~noncommutative extension of the problem raised by C.~Shannon since he was concerned with scalar
valued functions and we are dealing with matrix valued ones.

The {\em time limiting operator}, at level~$N$ acts on $\ell^2(M_R,{\mathbb N}_0)$ by simply setting equal to zero all
the components with index larger than~$N$.
We denote it by $\chi_N$.
The {\em band limiting operator}, at level~$\alpha$, acts on $L^2(W)$ by multiplication by the characteristic function
of the interval $(a, \alpha)$, $\alpha\le b$.
This operator will be denoted by $\chi_\alpha$.

One could consider restricting the band to an arbitrary subinterval $(a_1,b_1)$.
However, the algebraic properties exhibited here, see Section~\ref{Sec-commM} and beyond, hold only with this restriction.
A~similar situation arises in the classical case going all the way back to Shannon.

Consider the problem of determining a~function~$f$, from the following data:~$f$ has support on the f\/inite set
$\{0,\dots, N\}$ and its Fourier transform $\mathcal Ff$ is known on a~compact set $[a,\alpha]$.
This can be formalized as follows
\begin{gather*}
\chi_\alpha \mathcal Ff=g=\text{known},
\qquad
\chi_N f=f.
\end{gather*}
We can combine the two equations into
\begin{gather*}
E f= \chi_\alpha \mathcal F \chi_N f=g.
\end{gather*}
To analyze this problem we need to compute the singular vectors (and values) of the operator $E\colon \ell^2(M_R,{\mathbb
N}_0)\longrightarrow L^2(W) $.
These are given by the eigenvectors of the operators
\begin{gather*}
E^*E= \chi_N \mathcal F^{-1} \chi_\alpha \mathcal F \chi_N
\qquad
\text{and}
\qquad
S_2=E E^*= \chi_\alpha \mathcal F \chi_N \mathcal F^{-1} \chi_\alpha.
\end{gather*}
The operator $E^*E$, acting in $\ell^2(M_R,{\mathbb N}_0)$ is just a~f\/inite dimensional block-matrix~$M$, and each block
is given~by
\begin{gather*}
(M)_{m,n}=(E^*E)_{m,n}= \int_a^\alpha Q_m(x) W(x) Q{^*}_nw(x) dx,
\qquad
0\leq m,n \leq N.
\end{gather*}
The second operator $S_2= E E^*$ acts in $L^2((a,\alpha), W(t)dt)$ by means of the integral kernel
\begin{gather*}
k(x,y)=\sum\limits_{w=0}^N Q_w(x)Q{^*}_w(y).
\end{gather*}
Consider now the problem of f\/inding the eigenfunctions of $E^*E$ and $E E^*$.
For arbitrary~$N$ and~$\alpha$ there is no hope of doing this analytically, and one has to resort to numerical methods
and this is not an easy problem.
Of all the strategies one can dream of solving this problem, none sounds so appealing as that of f\/inding an operator
with simple spectrum which would have the same eigenfunctions as the original operators.
This is exactly what Slepian, Landau and Pollak did in the scalar case, when dealing with the real line and the actual
Fourier transform.
They discovered (the analog of) the following properties:
\begin{itemize}\itemsep=0pt
\item For each~$N$,~$\alpha$ there exists a~symmetric tridiagonal matrix~$L$, with simple spectrum, commuting with~$M$.
\item For each~$N$,~$\alpha$ there exists a~selfadjoint dif\/ferential operator~$D$, with simple spectrum, commuting with
the integral operator $S_2=EE^*$.
\end{itemize}

\section{From the real line to the sphere for matrix valued functions}

In this paper the role of the real line will be taken up by the~$n$-dimensional sphere.
We will consider $2 \times 2$ matrix valued functions def\/ined on the sphere with the appropriate invariance that makes
them functions of the colatitude~$\theta$ and we will use $x=\cos(\theta)$ as the variable.
The role of the Fourier transform will be taken by the expansion of our functions in terms of a~basis of matrix valued
orthogonal polynomials described below.
This is similar to the situation discussed in~\cite{GLP} except for the crucial fact that our functions are now matrix valued.
This situation has, to the best of our knowledge, not been considered before.

The matrix valued orthogonal polynomials considered here are those studied in~\cite{PZ}, arising from the spherical
functions of fundamental representations associated to the~$n$-dimensional sphere $S^n\simeq G/K$, where
$(G,K)=({\mathrm{SO}}(n+1), {\mathrm{SO}}(n))$, studied in~\cite{TZ}.

These spherical functions give rise to sequences $\{P_w\}_{w\geq 0}$ of matrix orthogonal polynomials depending on two
parameters~$n$ and $p \in \mathbb R$ such that $0< p< n$
\begin{gather*}
%\label{Pwdef}
P_w(x)=
\begin{pmatrix}
\dfrac{1}{n+1} C_w^{\frac{n+1}{2}}(x)+\frac{1}{p+w} C_{w-2}^{\frac{n+3}{2}}(x)&\dfrac{1}{p+w} C_{w-1}^{\frac{n+3}{2}}(x)
\vspace{2mm}\\
\dfrac{1}{n-p+w} C_{w-1}^{\frac{n+3}{2}}(x)&\dfrac{1}{n+1} C_w^{\frac{n+1}{2}}(x)+\dfrac{1}{n-p+w}
C_{w-2}^{\frac{n+3}{2}}(x)
\end{pmatrix}
,
\end{gather*}
where $C_w^\lambda(x)$ denotes the~$w$-th Gegenbauer polynomial
\begin{gather*}
C_w^\lambda(x)=\frac{(2\lambda)_w} {w!} \, {}_2 F_1\left(
\begin{matrix}
-w, w+2\lambda
\\
\lambda+1/2
\end{matrix}
;\frac{1-x}{2}\right),
\qquad
x\in[-1,1].
\end{gather*}
We recall that $C_w^\lambda$ is a~polynomial of degree~$w$ whose leading coef\/f\/icient is $\frac{2^w(\lambda)_w}{w!}$,
where $(a)_w=a(a+1)\cdots (a+w-1)$ denotes the Pochhammer's symbol.

In particular we have
\begin{gather*}
P_0=\frac 1{n+1}I,
\qquad
P_1=
\begin{pmatrix}
x & \frac 1{p+1}
\vspace{1mm}\\
\dfrac 1{n-p+1} & x
\end{pmatrix}
,
\\
P_2=
\begin{pmatrix}
\dfrac{(n+3)}2 x^2-\dfrac p{2(p+2)} & \dfrac{(n+3)}{p+2}x
\vspace{2mm}\\
\dfrac{(n+3)}{n-p+2}x & \dfrac{(n+3)}2 x^2-\dfrac {n-p}{2(n-p+2)}
\end{pmatrix}
.
\end{gather*}

Let us observe that the $\deg(P_w)=w$ and the leading coef\/f\/icient of $P_w$ is a~non singular scalar matrix
\begin{gather*}%\label{leadcoef}
\frac {2^w \big(\tfrac{n+1}2\big)_w}{(n+1) w!} \text{Id}.
\end{gather*}
The matrix polynomials $\{P_w\}_{w\geq 0} $ are orthogonal with respect to the matrix valued inner product
\begin{gather*}
\langle P,Q \rangle = \int_{-1}^1 P(x) W(x) Q(x)^* dx,
\end{gather*}
where the weight matrix is given by
\begin{gather}
\label{peso-x}
W(x)=W_{p,n}(x)= \big(1-x^2\big)^{\tfrac n2 -1}
\begin{pmatrix}
p x^2+n-p & -nx
\\
-nx & (n-p)x^2+p
\end{pmatrix}
,
\qquad
x\in [-1,1].
\end{gather}

Let us observe that by changing~$p$ by $n-p$, the weight matrices are conjugated to each other.
In fact, by taking $J=\left(
\begin{smallmatrix}
0&1
\\
1&0
\end{smallmatrix}
\right)$ we get
\begin{gather}
\label{p/n-ppeso}
J W_{p,n}J^* =W_{{n-p},n}.
\end{gather}

{\sloppy As a~consequence of this fact (or directly from the explicit def\/inition of $P_w$) we have that~$(P_w)_{22}$, the entry
$(2,2)$ of~$P_w$, is the same that the entry $(P_w)_{11}$ by replacing~$p$ by $n-p$.
Also the entry~$(2,1)$ of~$P_w$, $(P_w)_{21}$ is the entry $(P_w)_{12}$, if we replace~$p$ by $n-p$.

}

We have that $\langle P_w,P_w\rangle$ is always a~diagonal matrix.
Moreover one can verify that
\begin{gather*}
\langle P_w,P_w\rangle   =\| P_w\|^2= \frac{\sqrt \pi 2^{[w/2]}\Gamma\big(\tfrac n2+1+[\tfrac w2]\big) \prod\limits_{k=0}^{[(w-1)/2]}(n+2k+1)}{w!
(n+1)(n+2w+1)\Gamma\big(\tfrac n2+\tfrac 32\big)}
\\
\hphantom{\langle P_w,P_w\rangle   =\| P_w\|^2= }{}\times
\begin{pmatrix}
\dfrac {p (n-p+w+1)}{p+w}&0
\\
0& \dfrac{(n-p)(p+w+1)}{n-p+w}
\end{pmatrix}
.
\end{gather*}

We consider the orthonormal sequence of matrix polynomials
\begin{gather}
\label{relationQP}
Q_w= S_wP_w,
\end{gather}
where $S_w=\|P_w\|^{-1}$ is the inverse of the matrix $\|P_w\|$.

We display the f\/irst elements of the sequence $\{Q_w\}$.
\begin{gather}
Q_0  =\sqrt{\dfrac{\Gamma\big(\frac{n}{2}+ \frac{3}{2}\big)}{\sqrt \pi \Gamma\big(\frac{n}{2}+1\big)}}
\begin{pmatrix}
\dfrac 1{\sqrt{n-p+1}} &0
\\
0& \dfrac 1{\sqrt{p+1}}
\end{pmatrix},
\nonumber
\\
Q_1 = \sqrt {{\dfrac {2 \Gamma \big(\frac{n}{2}+\frac{5}{2} \big)}{\sqrt{\pi} \Gamma \big(\frac{n}{2}+1 \big)}}}
\begin{pmatrix}
\sqrt {{\dfrac {(p+1)}{p (n-p+2)}}} x & \dfrac 1 {\sqrt {p (n-p+2)(p+1)}}
\vspace{1mm}\\
\dfrac 1{\sqrt {(n-p) (p+2)(n-p+1)}} & \sqrt {\dfrac {(n-p+1)}{(n-p) (p+2)}} x
\end{pmatrix}.
\label{Q0}
\end{gather}

\section[The matrix~$M$]{The matrix~$\boldsymbol{M}$}%\label{EME}

Given the sequence of matrix orthogonal polynomials $\{P_w\}_{w\geq 0}$ we f\/ix a~natural number~$N$ and $\alpha \in
(-1,1)$ and consider the matrix~$M$ of total size $2(N+1)\times 2(N+1)$,
\begin{gather*}
M=
\begin{pmatrix}
M^{0,0}& M^{0,1}& \dots & M^{0,N}
\\
M^{1,0}& M^{1,1}& \dots & M^{1,N}
\\
\dots & \dots & \dots & \dots
\\
M^{N,0}& M^{N,1}& \dots & M^{N,N}
\end{pmatrix}
,
\end{gather*}
whose $(i,j)$-block is the $2\times 2$ matrix obtained by taking the inner product of the~$i$-th and~$j$-th {\em
normalized} matrix valued orthogonal polynomials in the interval $[-1,\alpha]$ with $\alpha\leq 1$
\begin{gather}
\label{MwithP}
M^{i,j}=\int_{-1}^\alpha Q_i(x) W(x) Q_j(x)^* dx
\qquad
\text{for}
\quad
0\leq i,j\leq N.
\end{gather}
It should be clear that the restriction to the interval $[-1,\alpha]$ implements ``band-limiting'' while the restriction
to the range $0, 1,\dots, N$ takes care of ``time-limiting''.
In the language of~\cite{GLP} where we were dealing with scalar valued functions def\/ined on spheres the f\/irst
restriction gives us a~``spherical cap'' while the second one amounts to truncating the expansion in spherical harmonics.

We gather here a~few important properties of the matrix $M$.

The entries of the matrix $M=(M_{rs})_{1\leq r,s\leq 2(N+1)}$ are related with the entries of the block matrices
$M^{i,j}=\left(
\begin{smallmatrix}
M^{i,j}_{11}& M^{i,j}_{12}
\\
M^{i,j}_{21}& M^{i,j}_{22}
\end{smallmatrix}
\right)$~by
\begin{gather*}%\label{blocks/entries}
M_{2i+1,2j+1} = M^{i,j}_{11},
\qquad\!
M_{2i+1,2j}= M^{i,j-1}_{12},
\qquad\!
M_{2i,2j+1} = M^{i-1,j}_{21},
\qquad\!
M_{2i,2j} = M^{i-1,j-1}_{22}.
\end{gather*}

From the def\/inition~\eqref{MwithP} it is clear that since $M^{j,i}=(M^{i,j})^*$,~$M$ is a~symmetric matrix
\begin{gather*}
M=M^*.
\end{gather*}

The weight matrices $W_{p,n}$ and $W_{n-p,p}$ are conjugated to each other by~\eqref{p/n-ppeso}.
Let us denote $M^{i,j}(p)$ the $2\times 2$ matrix with parameter~$p$, we get that the entry $1$, $2$ of $M^{i,j}$ is equal
to the entry $2$, $1$ by replacing~$p$ by $n-p$, i.e.,
\begin{gather}
\label{relationinblocks}
M^{i,j}_{12}(p)= M^{i,j}_{21}(n-p),
\qquad
M^{i,j}_{22}(p)= M^{i,j}_{11}(n-p).
\end{gather}
In Section~\ref{exp} we will give detailed information about the entries of~$M$.

\section[The commutant of~$M$ contains block tridiagonal matrices]{The commutant of~$\boldsymbol{M}$ contains block tridiagonal matrices}
\label{Sec-commM}

The aim of this section is to f\/ind all block tridiagonal symmetric matrices~$L$ such that
\begin{gather*}
L M= ML.
\end{gather*}

Notice that in principle there is no guarantee that we will f\/ind any such~$L$ except for a~scalar multiple of the identity.
For the problem at hand we need to exhibit matrices~$L$ that have a~simple spectrum, so that its eigenvectors will
automatically be eigenvectors of $M$.

Our f\/inding is that the vector space of such matrices~$L$ is of dimension 4.
It consist of the linear span of the matrices $L_1$, $L_2$ and $L_3$ given below.
Of course we can always add a~multiple of the identity matrix.
We do not have at the moment nice and clean proofs of these facts.
What we do have are careful symbolic computations, done (independently) in Maxima and in Maple, showing that the space
of solutions of the equation $ML=LM$ is always of dimension four and spanned by $L_1$, $L_2$, $L_3$ and the identity matrix.

We are only indicating explicitly the non-zero entries of these symmetric matrices, for $k=1, \dots,N+1$ we have
\begin{gather*}
%\label{eles}
(L_1)_{{2k,2k}} =\frac {(n-p+k-1)(p+k)}{(p+1)(n-p)},\displaybreak[0]
\\
(L_1)_{{2k-1,2k}}   =\frac{- (n-p+N+1)(p+N+1)}{\alpha (p+1)(n-2p)}\sqrt {\frac {p (p+k)(n-p+k-1)}{(n-p+k)(p+k-1)(n-p)}},
\displaybreak[0]
\\
(L_1)_{{2k,2k+2}} =\frac{(N-k+1) (N+n+k+1)}{{\alpha} (p+1)(n-p)} \sqrt{{\frac {k
(n-p+k-1)(p+1+k)(n+k)}{(n+2k-1)(n+2k+1)(p+k)(n-p+k)}}},
\\
(L_2)_{2k-1,2k-1} =\frac{(n-p+k)(p+k-1)}{p (n-p+1)},
\\
{(L_2)}_{{2k-1,2k}}  = \frac{(n-p+N+1)(p+N+1)}{\alpha (n-p+1)(n-2p)}\sqrt {\frac {(n-p+k)(p+k-1)(n-p)}{p
(p+k)(n-p+k-1)}},
\\
{(L_2)}_{{2k-1,2k+1}} =\frac{(N-k+1)(N+n+k+1)}{\alpha p (n-p+1)} \\
\hphantom{{(L_2)}_{{2k-1,2k+1}} =}{}\times\sqrt {\frac {k (p+k-1)(n-p+1+k)(n+k)}{(n+2k-1)(n+2k+1)
(p+k)(n-p+k)}}.
\end{gather*}
and f\/inally
\begin{gather*}
(L_3)_{{2k,2k}} =\frac {(k-1) (n+k) p (n-p+1)}{(n+2)(p+1)(n-p)},
\\
(L_3)_{{2k-1,2k-1}}  ={\frac {(k-1) (n+k)}{n+2}},
\\
(L_3)_{{2k-1,2k}}   = \frac{(k-1)(n-p+N+1)(p+N+1)(n+k) \sqrt p} {\alpha (p+1)(n+2)
\sqrt{(n-p)(p+k-1)(p+k)(n-p+k)(n-p+k-1)}},
\\
(L_3)_{{2k-1,2k+1}}  = \frac{(N-k+1)(N+n+k+1)}{\alpha (n+2)}   \frac {\sqrt{k
(p\!+\!k\!-\!1)(n\!-\!p\!+\!1\!+\!k)(n+k)}}{\sqrt{(n\!+\!2k\!-\!1)(n\!+\!2k\!+\!1)(p\!+\!k)(n\!-\!p\!+\!k)}},
\\
(L_3)_{{2k,2k+2}} = \frac{p (n-p+1)(N-k+1)(N+n+k+1){\sqrt{k (n\!-\!p\!+\!k\!-\!1)(p\!+\!k\!+\!1)(n\!+\!k)}}} {\alpha (n+2)(n-p)(p+1)
{\sqrt{ (n\!+\!2k\!-\!1  )  (n\!+\!2k\!+\!1  ) (p\!+\!k  )  (n\!-\!p\!+\!k  )}}}.
\end{gather*}

In terms of the $2\times 2$ blocks we have that
\begin{gather*}
L=
\begin{pmatrix}
L^{0,0} & L^{0,1}& 0 & \dots & 0
\\
L^{1,0}& L^{1,1} & L^{1,2} & \dots & 0
\\
0 & L^{2,1}& L^{2,2} &&
\\
\vdots && &\ddots & \vdots
\\
0 & \dots & && L^{N,N}
\end{pmatrix}
,
\end{gather*}
where the blocks are of the form
\begin{gather*}
\big(L_1\big)^{j,j} =
\begin{pmatrix}
0 & (L_1)_{2j+1,2j+2}
\\
(L_1)_{2j+1,2j+2} & (L_1)_{2j+2,2j+2}
\end{pmatrix}, \qquad
j=0, \dots, N,
\\
\big(L_1\big)^{j-1,j} =
\begin{pmatrix}
0 & 0
\\
0 & (L_1)_{2j,2(j+1)}
\end{pmatrix}, \qquad
j=1, \dots, N,
\\
\big(L_2\big)^{j,j} =
\begin{pmatrix}
(L_2)_{2j+1,2j+1} & (L_2)_{2j+1,2j+2}
\\
(L_2)_{2j+1,2j+2} & 0
\end{pmatrix}, \qquad
j=0, \dots, N,
\\
\big(L_2\big)^{j-1,j} =
\begin{pmatrix}
(L_2)_{2j-1,2j+1} & 0
\\
0 & 0
\end{pmatrix}, \qquad
j=1, \dots, N,
\\
\big(L_3\big)^{j,j} =
\begin{pmatrix}
(L_3)_{2j+1,2j+1} & (L_3)_{2j+1,2j+2}
\\
(L_3)_{2j+1,2j+2} & (L_3)_{2j+2,2j+2}
\end{pmatrix}, \qquad
j=0, \dots, N,
\\
\big(L_3\big)^{j-1,j} =
\begin{pmatrix}
(L_3)_{2j-1,2j+1} & 0
\\
0 & (L_3)_{2j,2(j+1)}
\end{pmatrix}, \qquad
j=1, \dots, N.
\end{gather*}

These matrices~$L$ are pentadiagonal matrices, moreover the of\/f-diagonal two-by-two blocks are diagonal matrices.

In the scalar case the matrix~$L$ is unique up to shifts and scaling, in the $2 \times 2$ case at hand this is no longer true.

\section[The spectrum of $L_1$]{The spectrum of $\boldsymbol{L_1}$}%\label{spectrum}

For the purpose at hand it is enough to exhibit one matrix with simple spectrum satisfying $ML=LM$.
Here we give a~``hands on'' argument to prove that the matrix $L_1$ has simple spectrum.
From this, and the commutativity established earlier, it follows that every eigenvector of $L_1$ is automatically an
eigenvector of~$M$, and we have achieved our goal: for each value of the parameters $(\alpha,N)$ we have a~numerically
ef\/f\/icient way to compute the eigenvectors of~$M$.
As a~referee has suggested, a~similar argument could be made to work in the case of $L_2$, $L_3$.

First observe the structure of the matrix $L_1$.
For $N=3$ we have
\begin{gather*}
L_1=\left(
\begin{array}
{cc|cc|cc|cc} 0& \ell_{12} & 0 &0 & 0 & 0 &0&0
\\
\ell_{12}& \ell_{22} & 0& \ell_{24} &0 & 0 &0&0
\\
\hline
0& 0 & 0 &\ell_{34} & 0 & 0&0&0
\\
0& \ell_{24} & \ell_{34} &\ell_{44} & 0 & \ell_{46}&0&0
\\
\hline
0& 0 & 0 &0 & 0 & \ell_{56}&0&0
\\
0& 0 & 0 & \ell_{46} & \ell_{56} & \ell_{66} &0&\ell_{68}
\\
\hline
0& 0 & 0 &0 & 0 & 0&0&\ell_{78}
\\
0& 0 & 0 &0 & 0 & \ell_{68} &\ell_{78}&\ell_{88}
\end{array}
\right),
\end{gather*}
where all indicated entries $\ell_{i,j}=(L_1)_{i,j}$ are nonzero,
see Section~\ref{Sec-commM}.

For general~$N$, with the notation $E_{i,j}$ for the standard basis of matrices, we have that $L_1$ is of the form
\begin{gather*}
L_1=\sum\limits_{j=1}^{N+1} \ell_{2j,2j} E_{2j,2j} +\ell_{2j-1,2j} \left(E_{2j-1,2j}+ E_{2j,2j-1}\right)
+\ell_{2j,2j+2}\left(E_{2j,2j+2}+ E_{2j+2,2j}\right)
\end{gather*}
(note that $\ell_{0,2}=\ell_{2N+2,2N+4}=0$).
Its $2\times 2$ blocks are of the form
\begin{gather*}
(L_1)^{j,j} =
\begin{pmatrix}
0 & \ell_{2j+1,2j+2}
\\
\ell_{2j+1,2j+2} & \ell_{2j+2,2j+2}
\end{pmatrix}
,
\qquad
(L_1)^{j-1,j} =
\begin{pmatrix}
0 & 0
\\
0 & \ell_{2j,2j+2}
\end{pmatrix}
,
\end{gather*}
where the coef\/f\/icients $\ell_{i,j}=(L_1)_{i,j}$ were given in Section~\ref{Sec-commM}.

One can see by induction that
\begin{gather*}
%\label{det}
\det (L_1) = (-1)^{N+1}\big (\ell_{1,2} \ell_{3,4} \cdots \ell_{2N+1,2N+2}\big)^2=
(-1)^{N+1}\prod\limits_{j=1}^{N+1}\big(\ell_{2j-1, 2j}\big)^2.
\end{gather*}

Assume that~$\lambda$ is an eigenvalue of $L_1$ and thus non-zero, and denote by~$X$ one of its eigenvectors.
We will show that~$X$ is a~scalar multiple of a~certain vector that depends only on the matrix $L_1$ and the eigenvalue
in question.
This shows that the geometric multiplicity of~$\lambda$ is one.

If $(x_1,x_2,x_3,\dots,x_{2N+2})$ are the components of~$X$, the scalar equations resulting from
\begin{gather*}
L_1 X = \lambda X
\end{gather*}
break up into two dif\/ferent groups: the equations given by the odd entries, $1,3,\dots,2N+1$, of the vector equation
above, namely
\begin{gather}
\label{oddeq}
\ell_{2j-1,2j} x_{2j} -\lambda x_{2j-1} =0,
\qquad
j=1, \dots,N+1,
\end{gather}
and the equations given by the even entries, $2,4,\dots,2N+2$, namely{\samepage
\begin{gather*}
%\label{eveneq}
\ell_{2j,2j+2} x_{2j+2}+ (\ell_{2j,2j} -\lambda) x_{2j} + \ell_{2j-1,2j} x_{2j-1} + \ell_{2j-2,2j} x_{2j-2} = 0,
\end{gather*}
for $j= 1,\dots,N+1$, with the convention $\ell_{0,2}=\ell_{2N+2,2N+4}=0$.}

The f\/irst set of equations allows us to write $x_1,x_3,x_5,\dots,x_{2N+1}$ in terms of $x_2,x_4,x_6,\dots$, $x_{2N+2}$.
We get
\begin{gather}
\label{xodd}
x_{2j-1} = \frac {\ell_{2j-1,2j}}{\lambda} x_{2j},
\qquad
j=1,2,\dots,N+1.
\end{gather}

By replacing~\eqref{xodd} in~\eqref{oddeq} we get, for $j=1, \dots, N+1$,
\begin{gather}
\label{eveneq2}
\ell_{2j,2j+2} x_{2j+2}+ \left(\ell_{2j,2j} -\lambda + \dfrac {(\ell_{2j-1,2j})^2}{\lambda} \right) x_{2j} + \ell_{2j-2,2j}
x_{2j-2} = 0.
\end{gather}

This set of equations allows us to f\/ind, successively, the entries $x_{2j+2}$ in terms of $x_2$, for $j=1,\dots,N$.
For $j=1$ we obtain
\begin{gather*}
x_4=-\frac 1{\ell_{2,4}}\left(\ell_{2,2}-\lambda+\frac{(\ell_{1,2})^2}\lambda\right) x_2.
\end{gather*}

For $j=2$ and using this expression for $x_4$, we get
\begin{gather*}
x_6=\frac 1{\ell_{2,4}\ell_{4,6}} \left(-(\ell_{2,4})^2 +\left(\ell_{4,4}-\lambda+\frac{(\ell_{3,4})^2}\lambda
\right)\left(\ell_{2,2}-\lambda+\frac{(\ell_{1,2})^2}\lambda\right)\right) x_2.
\end{gather*}

Our strategy is now clear, by successively solving the equations~\eqref{eveneq2} and using the relation between~$x_{2j}$
and~$x_{2j-1}$ we obtain expressions for~$x_2, x_3, x_4,\dots$ etc.
as quantities built out of~$L_1$ and~$\lambda$ all multiplied by the free parameter~$x_1$.

This goes on till we try to solve the last equation in~\eqref{eveneq2}, with $j=N+1$, where we meet our only
restriction, namely
\begin{gather*}
\left(\ell_{2N+2,2N+2} -\lambda + \frac {(\ell_{2N+1,2N+2})^2}{\lambda} \right) x_{2N+1} + \ell_{2N,2N+2} x_{2N} = 0.
\end{gather*}

With the variables $x_2, x_3, x_4,\dots,x_{2N+1}, x_{2N+2}$ all expressed in the desired form (i.e., as certain multiples of~$x_1$)
the last equation becomes a~polynomial in elements of~$L_1$ and the unknown quantity~$\lambda$ all multiplied by~$x_1$.
Except for this factor this is just the characteristic polynomial for $L_1$ and by assumption~$\lambda$ is a~root of it.

\section[Explicit expression for the matrix~$M$]{Explicit expression for the matrix~$\boldsymbol{M}$}\label{exp}

Writing down nice and explicit expressions for the entries of~$M$ is not easy.
We display here some of the entries whose expressions are not too involved.

We start with computing explicitly the entries
of the matrix $M^{0,0}=\displaystyle \int_{-1}^\alpha Q_0(x) W(x) Q_0(x)^*dx$.
By the def\/inition of the weight matrix given in~\eqref{peso-x} and the explicit expression of the $Q_0$ given
in~\eqref{Q0} we have
\begin{gather*}
M^{0,0} = \frac{\Gamma\big(\frac n2+\frac 32\big)}{\sqrt \pi \Gamma\big(\frac n2+1\big)} \int_{-1}^\alpha \big(1-x^2\big)^{\frac n2-1}
\begin{pmatrix}
\dfrac{px^2+n-p}{n-p+1}& -\dfrac{nx}{\sqrt{(p+1)(n-p+1)}}
\vspace{1mm}\\
-\dfrac{nx}{\sqrt{(p+1)(n-p+1)}} & \dfrac{(n-p)x^2+p}{p+1}
\end{pmatrix}
\\
\hphantom{M^{0,0}}{}
=
\begin{pmatrix}
M^{0,0}_{11}& M^{0,0}_{12}
\\
M^{0,0}_{21} &M^{0,0}_{22}
\end{pmatrix}
.
\end{gather*}

It is easy to verify that
\begin{gather*}
M^{0,0}_{12}= M^{0,0}_{21}=\frac{\Gamma\big(\frac n2+\frac 32\big) \big(1-\alpha^2\big)^{\frac n2}}{\sqrt {\pi (p+1)(n-p+1)}
\Gamma\big(\frac n2+1\big)}.
\end{gather*}

We also have
\begin{gather*}
M^{0,0}_{11} = \frac{\Gamma\big(\frac n2+\frac 32\big)}{(n-p+1) \sqrt \pi \Gamma\big(\frac n2+1\big)}\\
\left.\hphantom{M^{0,0}_{11} =}{}\times
\left(\tfrac 13 p x^3 \, {}_2 F_1\left(
\begin{smallmatrix}
-\frac n2 +1, \frac 32
\\
\frac 52
\end{smallmatrix}
;x^2\right) + (n-p) x \, {}_2 F_1\left(
\begin{smallmatrix}
-\frac n2 +1, \frac 12
\\
\frac 32
\end{smallmatrix}
;x^2\right)\right) \right|_{-1}^\alpha.
\end{gather*}
Computing the hypergeometric function at $x=1$ we get
\begin{gather}
\label{M11}
M^{0,0}_{11} = \frac{\Gamma\big(\frac n2+\frac 32\big) \Big(\tfrac 13 p \alpha^3 \, {}_2 F_1\left(
\begin{smallmatrix}
- n/2 +1, 3/2
\\
5/2
\end{smallmatrix}
;\alpha^2\right) + (n-p) \alpha\,  {}_2 F_1\left(
\begin{smallmatrix}
- n/2 +1, 1/2
\\
3/2
\end{smallmatrix}
;\alpha^2\right) \Big)}{(n-p+1) \sqrt \pi \Gamma\big(\frac n2+1\big)} + \frac 12.
\end{gather}

The expression of $M^{0,0}_{22}$ can be obtained from $M^{0,0}_{11}$ by changing~$p$ by $n-p$,
see~\eqref{relationinblocks}.

\subsection[The matrix $\widetilde M$]{The matrix $\boldsymbol{\widetilde M}$}

To simplify some calculations that follow and to avoid some square roots, we consider the following matrix $\widetilde M$.
This matrix is close to the matrix~$M$ and it is def\/ined as a~matrix of size $2(N+1)\times 2(N+1)$, whose $(i,j)$-block
is the $2\times 2$ matrix given~by
\begin{gather*}
\widetilde M^{i,j}=\int_{-1}^\alpha P_i(x) W(x) P_j(x)^* dx
\qquad
\text{for}
\quad
0\leq i,j\leq N.
\end{gather*}
We observe that $\widetilde M$ is def\/ined in the same way as~$M$ by using the sequence $(P_w)_{w\geq 0}$ of matrix
orthogonal polynomial instead of the orthonormal polynomials $(Q_w)_{w\geq 0}$.

The sequences of orthogonal polynomials $(P_w)_{w\geq 0}$ and $(Q_w)_{w\geq 0}$ are related by $Q_w=S_wP_w$, where
$S_w=\|P_w\|^{-1}$ (see~\eqref{relationQP}).
Thus the blocks of the matrices~$M$ and $\widetilde M$ are related~by
\begin{gather*}
M^{i,j}  =S_i \tilde M^{i,j} S_j^*
\qquad
\text{for}
\quad
0\leq i,j\leq N.
\end{gather*}

Since the matrices $S_j$ are diagonal matrices, we have that
\begin{gather*}
%\label{relationM-Mtilde}
M_{11}^{i,j} =(S_i)_{11}(S_j)_{11}\widetilde M^{i,j}_{11},
\qquad
M_{22}^{i,j}=(S_i)_{22}(S_j)_{22}\widetilde M^{i,j}_{22},
\\
M_{12}^{i,j} =(S_i)_{11}(S_j)_{22}\widetilde M^{i,j}_{12},
\qquad
M_{21}^{i,j}=(S_i)_{22}(S_j)_{11}\widetilde M^{i,j}_{21}.
\end{gather*}

Explicitly,
\begin{gather*}
S_j=\sqrt{\frac{j! (n+1)(n+2j+1)\Gamma\big(\tfrac n2+\tfrac 32\big)} {\sqrt \pi 2^{[j/2]}\Gamma\big(\tfrac n2+1+[\tfrac j2]\big)
\prod\limits_{k=0}^{[(j-1)/2]}(n+2k+1)}}\\
\hphantom{S_j=}{}\times
\begin{pmatrix}
\sqrt{\dfrac{p+j}{p (n-p+j+1)}} &0
\\
0& \sqrt{\dfrac{n-p+j}{(n-p)(p+j+1)}}
\end{pmatrix}
.
\end{gather*}

\subsection[Non diagonal elements for all the blocks of~$M$]{Non diagonal elements for all the blocks of~$\boldsymbol{M}$}

By direct computations we obtained the following expression for the non diagonal elements of each $2\times 2$ block $\widetilde M^{i,j}$.

For $0\leq i$, $j\leq N$ we get
\begin{gather*}
%\label{nondiagonalMentries}
\widetilde M_{1,2}^{i,j}= \frac{\big(1-\alpha^2\big)^{n/2} p (n-p)}{(n+1)^2(p+i)(n-p+j)} C_i^{\tfrac{n+1}2} (\alpha)
C_j^{\tfrac{n+1}2} (\alpha),
\\
\widetilde M_{2,1}^{i,j}=\frac{\big(1-\alpha^2\big)^{n/2} p (n-p)}{(n+1)^2(p+j)(n-p+i)} C_i^{\tfrac{n+1}2} (\alpha)
C_j^{\tfrac{n+1}2} (\alpha),
\end{gather*}
where $C_n^{\tfrac{n+1}2}(\alpha)$ denotes the~$n$-th Gegenbauer polynomial in~$\alpha$.

We observe that these expressions determine the entries $M_{2i+1,2j}$ and $M_{2i, 2j+1}$ of the mat\-rix~$M$,
for $1\leq i$, $j\leq N+1$, are given~by
\begin{gather*}
M_{2i+1,2j} = M_{1,2}^{i,j-1}= (S_i)_{11}(S_{j-1})_{22}\widetilde M_{1,2}^{i,j},
\qquad
M_{2i,2j+1} = M_{2,1}^{i-1,j}= (S_{i-1})_{22}(S_{j})_{11}\widetilde M_{2,1}^{i,j}.
\end{gather*}

\subsection[The diagonal entries of non diagonal blocks of~$M$]{The diagonal entries of non diagonal blocks of~$\boldsymbol{M}$}

Here we give the diagonal entries of the non diagonal blocks of $\widetilde M$, i.e., the elements $\widetilde M_{11}^{j,k}$
and~$\widetilde M_{22}^{j,k}$, with $j\neq k$.

For $1\leq j$, $k\leq N$, $j\neq k$ we have
\begin{gather*}
\widetilde M^{j-1,k-1}_{11}= \frac{p \big(1-\alpha^2\big)^{n/2}}{(k-j)(j+k+n-1)(n+1)^2} \left(\frac{j
(n-p+j-1)(n-p+k)}{n+2j-1} C_{j}(\alpha)C_{k-1}(\alpha) \right.
\\
\phantom{\widetilde M^{j-1,k-1}_{11}=}{}
- \frac{k (n-p+k-1)(n-p+j)}{n+2k-1} C_{j-1}(\alpha)C_{k}(\alpha)
\\
\phantom{\widetilde M^{j-1,k-1}_{11}=}{}
+\frac{(p+j)(n+j-1)(n-p+j)(n-p+k)}{(p+j-1)(n+2j-1)} C_{j-2}(\alpha)C_{k-1}(\alpha)
\\
\phantom{\widetilde M^{j-1,k-1}_{11}=}{}
\left.
-\frac{(p+k)(n+k-1)(n-p+j)(n-p+k)}{(p+k-1)(n+2k-1)} C_{j-1}(\alpha)C_{k-2}(\alpha) \right).
\end{gather*}
The element $\widetilde M_{22}^{j-1,k-1}$ is obtained from $\widetilde M_{11}^{j-1,k-1}$ by changing~$p$ by $n-p$.

These expressions determine the entries $M_{2j-1,2k-1}$ and $M_{2j, 2k}$ of the matrix~$M$, for $1\leq j\neq k \leq N+1$.
\begin{gather*}
M_{2j-1,2k-1} = M_{1,1}^{j-1,k-1}= (S_{j-1})_{11}(S_{k-1})_{11}\widetilde M_{1,1}^{j-1,k-1},
\\
M_{2j,2k} = M_{2,2}^{j-1,k-1}= (S_{j-1})_{22}(S_{k-1})_{22}\widetilde M_{2,2}^{j-1,k-1}.
\end{gather*}

\subsection[The main diagonal of $M$]{The main diagonal of $\boldsymbol{M}$}

At this point we have explicit expressions of all non diagonal entries of~$M$ (and $\widetilde M$).
Now we want to describe how to obtain the elements $M_{11}^{i,i}=M_{2i+1,2i+1}$ and $M_{22}^{i,i}= M_{2i+2, 2i+2}$, in
terms of the known coef\/f\/icients.

Starting with the equation $(ML-LM)_{2j,2j+2}=0$, with $L=L_1$ and $j=1,\dots, N$, we obtain
\begin{gather*}
M_{2j,2j}  - M_{2j-1,2j-1}
\\
\qquad
 = (L_{2j-1,2j})^{-1}\big(M_{2j-1,2j} L_{2j,2j}+ M_{2j-1,2j-2} L_{2j-2,2j}+ M_{2j-1,2j+2} L_{2j,2j+2} \big).
\end{gather*}

The equation $(ML-LM)_{2j,2j+2}=0$ is
\begin{gather*}
\begin{split}
&M_{2j+2,2j+2}  - M_{2j,2j} = (L_{2j,2j+2})^{-1}\big(M_{2j,2j+2} (L_{2j+2,2j+2}-L_{2j,2j}) + M_{2j,2j+1} L_{2j+1,2j+2}
\\
& \hphantom{M_{2j+2,2j+2}  - M_{2j,2j} =}{}
- M_{2j-1,2j+2} L_{2j-1,2j} + M_{2j,2j+4}L_{2j+2,2j+4}- M_{2j-2,2j+2}L_{2j-2,2j} \big).
\end{split}
\end{gather*}

By using~\eqref{relationinblocks} we have that the element $ M_{22}^{j-1,k-1}$ is obtained from $ M_{11}^{j-1,k-1}$~by
changing~$p$ into $n-p$.

Since we already know the explicit value of $M_{1,1}$ (see~\eqref{M11}), these expressions allow us to determine all the
diagonal elements of $M$.

\section{Some numerical results}

In this last section we display the results of some numerical computations.
This should make clear the importance of having found, as above, a~matrix such as $L_1$ for a~given~$M$.

Our point becomes very clear even if we use a~small value of~$N$ and a~value of~$\alpha$ pretty close to~$1$.
If we had chosen a~larger value of~$N$ the phenomenon in question would be present for an even larger range of values
of~$\alpha$.
For our illustration we choose $n=27$, $p=15$, and f\/inally $N=2$ and $\alpha=9/10$.

If our task is to compute the eigenvectors of~$M$ we can use the QR algorithm as implemented in LAPACK.
The results are recorded below, where we denote by~$X$ the matrix of eigenvectors (normalized and given as columns
of~$X$) and by~$D$ the diagonal matrix of eigenvalues.
The matrix~$X$ is
\begin{gather*}
X=\left(
\begin{matrix}
0.046636 & -.0318428& 0.294888 & -.953861& - 0.0995136 & -.0406174
\\
-.0424748 & -.0353152 & -.241903 & -.0756026&.852604 & -.429742
\\
.214609 & -.17084 &.649236 & -.216917&.0981185 &-.654981
\\
-.19526 & -.189172 & -.589859 & -.18561 &.494728 & 0.569823
\\
.706579 & -.64988 & -.215899& 0.0105047 & -.0312226&.174508
\\
-.642327 & -.714438 &.19614 &.0530856 & -.0883495 & -.171814
\end{matrix}
\right)
\end{gather*}

The question is: should we trust the result produced by this high quality numerical package?
\\
One could be quite satisf\/ied by observing that the dif\/ference
\begin{gather*}
M X- X D
\end{gather*}
is indeed very small.
On the other hand LAPACK reports for eigenvalues of~$M$, with appropriate rounding-of\/f, the values
$1.0$, $1.0$, $1.0$, $1.0$, $1.0$, $1.0$.
This should be a~red f\/lag.

Recall that the eigenvectors of $L_1$ should agree (up to order) with those of~$M$.
If we denote the matrix made up of the normalized eigenvectors of $L_1$ by~$Y$ we get
\begin{gather*}
Y=\left(
\begin{matrix}
.641473 &.227019 &.0318478 & -.674518 &.280734 & -.0466363
\\
.688247 &.247073 &.0350924 &.628845 & -.258435 &.0424743
\\
-.229364 &.613934 &.170768 &.280102 &.645602 & -.214607
\\
- 0.24584 &.667491 &.187976 & -.260873 & -.593721 &.195257
\\
0.028825 & -.172614 &.649883 & -.0405441 & -.214796 & -.706573
\\
.0308702 & -.187518 & 0.71478 &.0377295 &.197372&.642335
\end{matrix}
\right)
\end{gather*}

For the eigenvalues of $L_1$, LAPACK returns the values $6.46314$, $6.55601$, $6.63761$, $-5.61601$, $-5.54541$, $-5.4863$
a~reasonably spread out spectrum.

If we compute the matrix of inner products given~by
\begin{gather*}
Y^T X
\end{gather*}
we expect to have the identity matrix up to some permutation and possibly some signs due to the normalization of the
eigenvectors which are the columns of~$X$ and~$Y$.
In our case we get for the moduli of the entries of~$Y^T X$ the matrix
\begin{gather*}
\left(
\begin{matrix}
4.65e^{-7} & 1.71e^{-4} &.019 &.667 &.671 &.021
\\
2.83e^{-6} & 9.61e^{-4} &.013 & 0.234 &.2244 &.9602
\\
9.85e^{-6} &.9999 & 7.801e^{-4} & 2.385e^{-4} & 4.395e^{-7} & 8.81e^{-4}
\\
7.49e^{-7} & 1.61e^{-4} & 8.569e^{-4} & 0.707 &.707 &.278
\\
3.93e^{-6} & 7.89e^{-4} &.9997 &.0147 &.007 &.0115
\\
1.0 & 9.84e^{-6} & 3.89e^{-6} & 9.61e^{-7} & 1.46e^{-6} & 2.55e^{-6}
\end{matrix}
\right)
\end{gather*}

The reader will observe that some of the entries of this matrix are indeed very close to the theoretically correct
values, while others are terribly of\/f.
The reason is that for this choice of~$\alpha$ there are a~few eigenvalues of~$M$ that are just too close together.
This produces numerical instability in the computation of the corresponding eigenvectors.
On the other hand all the eigenvalues of~$L_1$ are nicely separated, and the corresponding eigenvectors can be trusted.

In summary, a~good way to obtain good numerical values for the eigenvectors of~$M$ is to forget about~$M$ altogether and
to compute numerically the eigenvectors of~$L_1$.
Not only we will then be dealing with a~very sparse matrix for which the QR algorithm works very fast (most of the work
is avoided) but the problem is numerically very well conditioned.
For a~discussion of the sensitivity of eigenvectors and their dependence on the separation of the corresponding
eigenvalues one can consult~\cite{De,St} as well as~\cite[p.~15 and p.~222]{Pa}.

As a~referee has pointed out, these numerical problems are not new in our situation involving matrix valued functions,
but already appear in the classical scalar case, and this phenomenon is well documented.
The important point is that even in the matrix valued case we can exhibit commuting tridiagonal matrices that play the
same role of the ``prolate spheroidal dif\/ferential operator'' in the scalar case.

\subsection*{Acknowledgements}

This research was supported in part by the Applied Mathematical
Sciences subprogram of the Of\/f\/ice of Energy Research, USDOE, under Contract DE-AC03-76SF00098, by AFOSR grant
FA95501210087 through a~subcontract to Carnegie Mellon University, by CONICET grant PIP 112-200801-01533, by SeCyT-UNC
and by the Oberwolfach Leibniz Fellows Program.

\pdfbookmark[1]{References}{ref}
\LastPageEnding


\begin{thebibliography}{99}
\footnotesize \itemsep=0pt

\bibitem{BK14}
Bonami A., Karoui A., Uniform approximation and explicit estimates for the
  prolate spheroidal wave functions, \href{http://arxiv.org/abs/1405.3676}{arXiv:1405.3676}.

\bibitem{De}
Demmel J.W., Applied numerical linear algebra, \href{http://dx.doi.org/10.1137/1.9781611971446}{Society for Industrial and
  Applied Mathematics (SIAM)}, Philadelphia, PA, 1997.

\bibitem{DG}
Duistermaat J.J., Gr{\"u}nbaum F.A., Dif\/ferential equations in the spectral
  parameter, \href{http://dx.doi.org/10.1007/BF01206937}{\textit{Comm. Math. Phys.}} \textbf{103} (1986), 177--240.

\bibitem{G3}
Gr{\"u}nbaum F.A., A new property of reproducing kernels for classical
  orthogonal polynomials, \href{http://dx.doi.org/10.1016/0022-247X(83)90123-3}{\textit{J.~Math. Anal. Appl.}} \textbf{95} (1983),
  491--500.

\bibitem{G5}
Gr{\"u}nbaum F.A., Some new explorations into the mystery of time and band
  limiting, \href{http://dx.doi.org/10.1016/0196-8858(92)90015-O}{\textit{Adv. in Appl. Math.}} \textbf{13} (1992), 328--349.

\bibitem{G4}
Gr{\"u}nbaum F.A., Band-time-band limiting integral operators and commuting
  dif\/ferential operators, \textit{St.~Petersburg Math.~J.} \textbf{8} (1997),
  93--96.

\bibitem{G6}
Gr\"unbaum F.A., The bispectral problem: an overview, in Special Functions
  2000: Current Perspective and Future Directions ({T}empe, {AZ}), \href{http://dx.doi.org/10.1007/978-94-010-0818-1_6}{\textit{NATO
  Sci. Ser.~II Math. Phys. Chem.}}, Vol.~30, Kluwer Acad. Publ., Dordrecht,
  2001, 129--140.

\bibitem{GLP}
Gr{\"u}nbaum F.A., Longhi L., Perlstadt M., Dif\/ferential operators commuting
  with f\/inite convolution integral operators: some nonabelian examples,
  \href{http://dx.doi.org/10.1137/0142067}{\textit{SIAM~J. Appl. Math.}} \textbf{42} (1982), 941--955.

\bibitem{GY}
Gr{\"u}nbaum F.A., Yakimov M., The prolate spheroidal phenomenon as a
  consequence of bispectrality, in Superintegrability in Classical and Quantum
  Systems, \textit{CRM Proc. Lecture Notes}, Vol.~37, Amer. Math. Soc.,
  Providence, RI, 2004, 301--312, \href{http://arxiv.org/abs/math-ph/0303041}{math-ph/0303041}.

\bibitem{JB}
Jahn K., Bokor N., Revisiting the concentration problem of vector f\/ields within
  a spherical cap: a commuting dif\/ferential operator solution,
  \href{http://dx.doi.org/10.1007/s00041-014-9324-7}{\textit{J.~Fourier Anal. Appl.}} \textbf{20} (2014), 421--451,
  \href{http://arxiv.org/abs/1302.5261}{arXiv:1302.5261}.

\bibitem{JKS}
Jamming P., Karoui A., Spektor S., The approximation of almost time and band
  limited functions by their expansion in some orthogonal polynomial bases,
  \href{http://arxiv.org/abs/1501.03655}{arXiv:1501.03655}.

\bibitem{SLP2}
Landau H.J., Pollak H.O., Prolate spheroidal wave functions, {F}ourier analysis
  and uncertainty.~{II}, \href{http://dx.doi.org/10.1002/j.1538-7305.1961.tb03977.x}{\textit{Bell System Tech.~J.}} \textbf{40} (1961),
  65--84.

\bibitem{SLP3}
Landau H.J., Pollak H.O., Prolate spheroidal wave functions, {F}ourier analysis
  and uncertainty. {III}.~{T}he dimension of the space of essentially time- and
  band-limited signals, \href{http://dx.doi.org/10.1002/j.1538-7305.1962.tb03279.x}{\textit{Bell System Tech.~J.}} \textbf{41} (1962),
  1295--1336.

\bibitem{M}
Melkman A.A., {$n$}-widths and optimal interpolation of time- and band-limited
  functions, in Optimal Estimation in Approximation Theory ({P}roc. {I}nternat.
  {S}ympos., {F}reudenstadt, 1976), Editors C.A.~Michelli, T.~Rivlin, Plenum,
  New York, 1977, 55--68.

\bibitem{ORX}
Osipov A., Rokhlin V., Xiao H., Prolate spheroidal wave functions of order
  zero. Mathematical tools for bandlimited approximation, \href{http://dx.doi.org/10.1007/978-1-4614-8259-8}{\textit{Applied
  Mathematical Sciences}}, Vol.~187, Springer, New York, 2013.

\bibitem{PZ}
Pacharoni I., Zurrian I., Matrix ultraspherical polynomials: the $2\times 2$
  fundamental cases, \textit{Constr. Approx.}, to appear, \href{http://arxiv.org/abs/1309.6902}{arXiv:1309.6902}.

\bibitem{Pa}
Parlett B.N., The symmetric eigenvalue problem, \textit{Prentice-Hall Series in
  Computational Mathematics}, Prentice-Hall, Inc., Englewood Clif\/fs, N.J., 1980.

\bibitem{P1}
Perlstadt M., Chopped orthogonal polynomial expansions~-- some discrete cases,
  \href{http://dx.doi.org/10.1137/0604012}{\textit{SIAM~J. Algebraic Discrete Methods}} \textbf{4} (1983), 94--100.

\bibitem{P2}
Perlstadt M., A property of orthogonal polynomial families with polynomial
  duals, \href{http://dx.doi.org/10.1137/0515081}{\textit{SIAM~J. Math. Anal.}} \textbf{15} (1984), 1043--1054.

\bibitem{PS}
Plattner A., Simons F.J., Spatiospectral concentration of vector f\/ields on a
  sphere, \href{http://dx.doi.org/10.1016/j.acha.2012.12.001}{\textit{Appl. Comput. Harmon. Anal.}} \textbf{36} (2014), 1--22,
  \href{http://arxiv.org/abs/1306.3201}{arXiv:1306.3201}.

\bibitem{S}
Shannon C.E., A mathematical theory of communication, \href{http://dx.doi.org/10.1002/j.1538-7305.1948.tb01338.x}{\textit{Bell System
  Tech.~J.}} \textbf{27} (1948), 379--423.

\bibitem{S+}
Shannon C.E., A mathematical theory of communication, \href{http://dx.doi.org/10.1002/j.1538-7305.1948.tb00917.x}{\textit{Bell System
  Tech.~J.}} \textbf{27} (1948), 623--656.

\bibitem{SD}
Simons F.J., Dahlen F.A., Spherical Slepian functions on the polar gap in
  geodesy, \href{http://dx.doi.org/10.1111/j.1365-246X.2006.03065.x}{\textit{Geophys.~J. Int.}} \textbf{166} (2006), 1039--1061,
  \href{http://arxiv.org/abs/math.ST/0603271}{math.ST/0603271}.

\bibitem{SDW}
Simons F.J., Dahlen F.A., Wieczorek M.A., Spatiospectral concentration on a
  sphere, \href{http://dx.doi.org/10.1137/S0036144504445765}{\textit{SIAM Rev.}} \textbf{48} (2006), 504--536,
  \href{http://arxiv.org/abs/math.CA/0408424}{math.CA/0408424}.

\bibitem{SLP4}
Slepian D., Prolate spheroidal wave functions, {F}ourier analysis and
  uncertainity. {IV}.~{E}xtensions to many dimensions; generalized prolate
  spheroidal functions, \href{http://dx.doi.org/10.1002/j.1538-7305.1964.tb01037.x}{\textit{Bell System Tech.~J.}} \textbf{43} (1964),
  3009--3057.

\bibitem{S2}
Slepian D., On bandwidth, \href{http://dx.doi.org/10.1109/PROC.1976.10110}{\textit{Proc. IEEE}} \textbf{64} (1976), 292--300.

\bibitem{SLP5}
Slepian D., Prolate spheroidal wave functions, {F}ourier analysis and
  uncertainity. {IV}.~The discrete case, \href{http://dx.doi.org/10.1002/j.1538-7305.1978.tb02104.x}{\textit{Bell System Tech.~J.}}
  \textbf{57} (1978), 1371--1430.

\bibitem{S1}
Slepian D., Some comments on {F}ourier analysis, uncertainty and modeling,
  \href{http://dx.doi.org/10.1137/1025078}{\textit{SIAM Rev.}} \textbf{25} (1983), 379--393.

\bibitem{SLP1}
Slepian D., Pollak H.O., Prolate spheroidal wave functions, {F}ourier analysis
  and uncertainty.~{I}, \href{http://dx.doi.org/10.1002/j.1538-7305.1961.tb03976.x}{\textit{Bell System Tech.~J.}} \textbf{40} (1961),
  43--63.

\bibitem{St}
Stewart G.W., Matrix algorithms. {V}ol.~{II}. Eigensystems, \href{http://dx.doi.org/10.1137/1.9780898718058}{Society for
  Industrial and Applied Mathematics (SIAM)}, Philadelphia, PA, 2001.

\bibitem{TZ}
Tirao J.A., Zurri{\'a}n I.N., Spherical functions of fundamental {$K$}-types
  associated with the {$n$}-dimensional sphere, \href{http://dx.doi.org/10.3842/SIGMA.2014.071}{\textit{SIGMA}} \textbf{10}
  (2014), 071, 41~pages, \href{http://arxiv.org/abs/1312.0909}{arXiv:1312.0909}.

\bibitem{W}
Weyl H., The theory of groups and quantum mechanics, Dutton, New York, 1931.

\end{thebibliography}
\end{document}